\documentclass{article}
\usepackage{hyperref}
\usepackage{amsfonts}
\usepackage{epigraph}

\title{Metaphysics in Mathematics}
\author{John Engelsted Hester}
\date{January 2017}

\begin{document}

\maketitle

\abstract{The purpose of this project is to outline various philosophies on the metaphysics of mathematics that have been prominent since the time of Cantor, highlighting some biographical aspects that have influenced these ideas as well.  The main topics will be the independent existence of mathematics and controversies over the nature of the mathematical infinite.  After describing these ideas, I will demonstrate some of the similarities that are common to these philosophies, showing that despite reaching different conclusions many are guided by the same driving idea.}

\section{Introduction}

\epigraph{It is clear that there is no classification of the Universe that is not arbitrary and full of conjectures.}{Jorge Luis Borges}

\hspace{4ex}The aim of this paper is to outline and describe various perspectives and metaphysical viewpoints of some prominent philosophers and mathematicians interested in logic and set theory from its inception in Cantor's time to today.  
I will include some historical notes that are relevant to the progression of an individuals thought or provide necessary context.
It is important to remark that these viewpoints are metaphysical in nature rather than mathematical, for instance a formalist may regard some of the work of a set theorist as legitimate mathematics in a formal system, while remaining agnostic towards the topic's existence as a meaningful object in a non-symbolic sense.
An example of this is Robinson's relationship with Kurt G{\"o}del, who wanted Abraham Robinson as his successor at the Institute for Advanced Study, despite G{\"o}del being considered a Platonist and Abraham Robinson being considered a formalist\ \cite[p. 481]{Kana}.

These ideas will be presented with a focus on various authors attitudes towards mathematical ideas that exist independently (or perhaps not) of human thought as abstract objects, specifically the existence of models of set theory and the existence of cardinals that are independent of the axioms of some set theory with the axiom of infinity, or perhaps some other dubious axiom omitted.
The paper will begin with the viewpoints of Georg Cantor, commonly regarded as the creator of set theory, and then cover the viewpoints of the contemporary intuitionism in contradistinction.
Next G{\"o}del's view will be presented along with the popular idea of mathematical formalism.

Ultimately I hope to find some of the common threads of the thinkers in the conclusion, and identify the artifacts of their viewpoints that remain in the spectrum of modern mathematical viewpoints.  I hope to bolster the hypothesis that the variety of perspectives towards mathematics today have not been created in a vacuum or by the philosophizing of a single person, but through a dialogue that has been unfolding for centuries and will certainly continue to do so, with fragments of even the staunchest of opponents coalescing in time.  

\section{Cantor the religious Platonist}
\hspace{4ex}George Cantor was a German mathematician who is most famous for his proof of the existence of multiple magnitudes of infinite numbers in his paper "On a Property of the Collection of all Real Algebraic Numbers"\footnote{\url{https://math.byu.edu/~grant/courses/cantor1874.pdf}} and for his formulation of the famous Continuum Hypothesis, which states there does not exist an infinite cardinality between that of the natural numbers and the cardinality of the continuum.  
His mathematical life can roughly be broken into two stages\ \cite[Ch. 6]{Cant}.
The first is a period of intense mathematical activity during which Cantor produced most of the mathematical works for which he is known today, and is associated with a philosophy of mathematics that exhibits the seeds of what became formalism. 
The second is a religious phase during which Cantor largely abandoned the pursuit of mathematics and instead favored philosophy and theology. 
It was during this second stage, incited by what amounted to his rejection by the mathematical community, that Cantor developed his famous religious perspective.  

Cantor always viewed his mathematical work as reflecting properties of the corporeal world; in his mind the cardinality of the continuum not only only represented the number of points on a line in Euclidean space, but the size of the collection of points, material or not, on an extended curve in real space.  
This view that mathematical ideas refer to real things that exist independently of man will be referred to Platonism from here on.  
To Cantor, the transfinite ordinals existed as naturally as the number five with its many occurrences in nature, perhaps in the number of trees on a hilltop.  
Just as the trees and the corresponding subcollections of trees exist in a natural way, the transfinite types arise from the concept of number and continuity in physical space for Cantor\ \cite[p. 171]{Tapp}.

During his first stage of philosophical thought and years as a working mathematician, Cantor held metaphysical views that while controversial at the time, might not be as surprising today.
Holding that mathematics was justified by its apparent consistency, Cantor believed that his exploration of the transfinite was legitimate as it seemed to lack any inconsistencies or contradictions.
While he still had faith in the immanent reality of the multitude of infinite numbers disclosed by his theory, Cantor made a point of distinguishing between his metaphysical claims and the main fare of mathematics that he proposed.  Here, he appeals to the firm foundation of rational numbers and points out that he argues for the transfinite in a way that parallels the justification of the rationals and complex numbeers.
As Cantor states, quoted in $\cite{Cant}$:

\begin{quote}
    In particular, in introducing new numbers, mathematics is only obliged to give definitions of them, by which such a definiteness, and circumstances permitting, such a relation to the older numbers are conferred upon them that in given cases they can definitely be distinguished from one another.  As soon as a number satisfies all these conditions, it can and must be regarded as existent and real in mathematics.  Here I perceive the reason why one has to regard the rational and complex numbers as being just as thoroughly existent as the finite positive integers.
\end{quote}

In this period, Cantor faced dogged opposition from some contemporary mathematicians such as Kronecker, and tried to sway them to acknowledge his work in this way.
Here mathematical existence is distinguished from "real" existence, with real mathematical existence referring to the idea of number as it is established in the mind via a construction or some other device.  
Cantor's view of mathematics being justified in its apparent consistency rather than purported truth in relation to extensional objects or spaces is similar to the viewpoints of some formalists.
However, it is his unreserved belief that the reality of mathematics in the mind must and does coincide with the idea of measure and number for extensional, real objects that established him as a Platonist in today's view.

Cantor's thought was immensely impacted by a series of mental breaks and professional mishaps in the mid 1880s that set in to motion what could be considered the second portion of his career.  During this second stage, Cantor further refined his philosophical ideas, largely with a bent that was influenced by his frustration at the mathematical community for what he interpreted as widespread rejection.
He had long hoped for a position at the University of G{\"o}ttingen or in Berlin, but events in his academic life at the hands of his enemies, both imagined and not, had made this impossible.
Viewing his transfinite numbers as an exploration of God, Cantor saw himself as delving into the incomprehensibility of the religious absolute through the indexing of transfinite species.  Rather than being legitimized by its apparent consistency, Cantor now considered his mathematics to be immanently legitimized by the existence of God.  
Cantor shared this viewpoint with his contemporary George Gutberlet\ \cite[Ch. 6]{Cant}:
\begin{quote}
But in the absolute mind the entire sequence is always in actual consciousness, without any possibility of increase in the knowledge or contemplation of a new member of the sequence.
\end{quote}

For Gutberlet and Cantor, the infinite knowledge of God represented an embodification of the infinite that both justified and encouraged the mathematical exploration of the transfinite.  
Further roused by endorsements from the Catholic church and neo-Thomists such as Cardinal Franzelin \cite{Cant}, the deeply religious Cantor continued to frame the infinite in a religious perspective.  
This appreciation that had always been lacking in the mathematical community is apparent in another quotation from \cite{Cant}.  Here Cantor is referring to his unfulfilled wish of a more prominent position at a German university:

\begin{quote}
But now I thank God, the all-wise and all-good, that He always denied me the fulfillment of this wish, for He thereby constrained me, through a deeper penetration into theology to serve Him and His Holy Roman Catholic Church better than I have been able with my exclusive preoccupation with mathematics. 
\end{quote}

\section{Kronecker the finitist}

\begin{quote}
God made the integers, all else is the work of man. -Kronecker
\end{quote}

\hspace{4ex}  The source of many of Cantor's woes in his academic career, or at least as he perceived it, was the often vitriolic relationship he had with the contemporary German mathematician Leopold Kronecker. 
Kronecker is often (incorrectly, as I shall later explain) labelled as the founder of the idea of intuitionism, namely that mathematics should be based on the intuitive counting of natural numbers.  
As might be suggested by the above quotation, he hoped for an arithmetization of mathematics, so that everything could be verifiable using computational methods.

He is often dismissed due to this viewpoint being largely incompatible with much of mathematics.  This was true even in his time druing the late 19th century, not to mention today\ \cite{Edwa}.
Sadly he is often attributed with many apocryphal quotations that perhaps started as jokes exaggerating his viewpoint, but ultimately falsely reflect his views.
However, it is important to note that he was very productive as a reseracher and had many valued ideas in algebra and number theory.
His hope for arithmetization of mathematics is well suited to a mathematical world that has seen only an increase in the use of electronic proof assistants over the years, as is insightfully pointed out by Harold Edwards.  

During the period of Cantor and Kronecker, and since antiquity, mathematicians and philosophers emphasized the distinction between the notions of potential and actual infinities. 
To illustrate this idea, potential infinity is found in the sequence of natural numbers $0,1,2,...$ while the collection $\{x\in\mathbb{Z}:x\geq 0\}$ would be considered an actual infinity.  
Up until Cantor's time, most felt that the notion of an actual infinity had little meaning, and considered the idea of an actual infinity an antinomy, as some consider it apparent that the mind is unable to rationally cope with what is truly infinite.  
One exception to this is the ancient Pythagoreans, who held in comparison a more modern view towards the existence of the actual infinite.  
This distinction between the potential and actual infinite is rooted in the thought of Aristotle, and in particular Kronecker inherited this view.  

It was during the mid 1880s that Kronecker began publicly expressing his finitistic views towards the philosophy of mathematics to his colleagues.  
It is often said that Kronecker did not practice the philosophy he preached, because he often used symbols that denoted mathematical constructs that would today be interpreted as having the touch of the actual infinite.
Instances of this include the number $\pi$, $e$, and the notion of the collection of all natural numbers.

This accusation is unfair, as he used such symbols as shorthand for the corresponding object that could be reasoned with in a finitistic manner.
For instance, $\pi$ can be reckoned with through geometric appeal to the circle and when speaking of the infinity of the natural numbers, he is referring to the fact that every natural number has a successor, not to any meaningful collection that contains them all.
In fact, as a result of this it is sometimes said that Kronecker did not believe some of his mathematical results because of the fact that some of them could be interpreted as using a completed infinity, when in fact he was using them as an abbreviation for a process.
Much of Kronecker's viewpoint can be encapsulated in this excerpt from his obituary\ \cite{Obit}:

\begin{quote}
Concerning the rigor of notions (he) imposes highest requirements and tries to squeeze everything that should have a right of citizenship in Mathematics into the crystal clear and edgy form of number theory.
\end{quote}

\section{Brouwer the intuitionist}

\hspace{4ex}One of Kronecker's philosophical successors was the famous topologist and analyst L.E.J. Brouwer.
Similar accusations of insincerity towards Brouwer's commitment to the finite are levelled towards him, but this is largely because his viewpoint towards the epistemological basis of mathematics evolved throughout his life and went through phases.
This means that his perspective towards what mathematics meant and justifies it fluctuated.
It is suggested that towards the end of his life, Brouwer did not even believe his famous fixed point theorem\ \cite{Brow}, which states that every continuous function from a closed disk to itself has a fixed point.
The reason for this is that it does not admit a constructive finite proof\ \cite{Brat}.  

Brouwer had the opinion that mathematics was a construction of the human mind, and during his most mature stage believed that in order to have a sound basis for mathematics, it must be the case that everything is based on finite constructions and operations.
Otherwise there is an opportunity for error to be introduced, due to the human inability to properly conceive of the infinite.
This viewpoint is known as intuitionism, the name being rooted in the idea of dealing exclusively with the most intuitive of abstract mathematical constructions, the natural numbers.  
Although Brouwer and Kronecker shared this heavy focus on the counting numbers, Kronecker's thought is distinct from Brouwer's in that Kronecker held that the natural numbers are intimately real as divine creation, in contrast to Brouwer's synthetic point of view.
For this reason it is inaccurate to label Kronecker as the forefather of the intuitionists.

The intuitionistic idea of acknowledging mathematics as nothing but a construction of the mind has several huge implications.  
Immediately, mathematics loses its privileged position as being a discipline that can make any claims of relevance to the world, meaning some sort of world that exists with priority ahead of our senses.  
Rather, mathematics can only be true so far as it is made clearly apparent as true to the individual in regards to a specific system. 
This could be paraphrased by saying that the stock of mathematical entities is a real thing, for each person, and for humanity, but not real apart from them.   

Perhaps the most revolutionary and broadly insightful consequence is that from the intuitionistic point of view, one can no longer defend the law of the excluded middle\ \cite{Intui}.
The law of the excluded middle is the statement that for any proposition $P$, one has $P\vee\neg P$.
Of course in a finite system, it is possible to check whether or not $P$ is true, so Brouwer's concern only appears in infinite cases.
This is because of the committed constructivism of the intuitionistic viewpoint: as mathematics is a construction of the human mind, anything held to be mathematically true must be a finite construction or verifiable with a finite method.  

This rejection foreshadows later developments in mathematics, in which things that are at first obvious and intuitive such as the law of excluded middle are cast into doubt, and their alternatives explored.  The emphasis on constructivism is something that flourishes with computational power, and the belief that a simple and constructive proof is preferable to a proof that is not is widespread today.  This is an example of an argument from one of Brouwer's lectures\ \cite{Brou} that intends to cast doubt on the excluded middle from the intuitionist perspective, lucidly arguing against it:

\begin{quote}
But now let us pass to infinite systems and ask for instance if there exists a natural number n such that in the decimal expansion of pi the nth, (n+1)th, ..., (n+8)th and (n+9)th digits form a sequence 0123456789. This question, relating as it does to a so far not judgeable assertion, can be answered neither affirmatively nor negatively. But then, from the intuitionist point of view, because outside human thought there are no mathematical truths, the assertion that in the decimal expansion of pi a sequence 0123456789 either does or does not occur is devoid of sense.
\end{quote}

\section{After intuitionism}

\hspace{4ex}There are many today that hold ideas similar to those of Brouwer today.  For instance, Paul Taylor in reference to his book \textit{Practical Foundations of Mathematics} refers to "The Great Set-Theoretic Swindle" as the idea that mathematics needs completed infinities as opposed to potential infinities for the bulk of mathematics.  Continuing in this vein, he casts doubt on the idea of a collection as necessarily having to correspond with that of a set, highlighting the fact that the axioms of set theory are rooted in basic intuitions about what we would consider to be a collection, without having any claim to reality in themselves.  From this perspective, it is meaningful to talk about a system as long as you concede that the object of study is the syntactical manipulation of symbols, rather than an encounter with some sort of semantic meaning in an extrasensory mathematical reality.

Part of the so-called swindle is the claim that in order to meaningfully talk about any sort of collection, one must be talking about a set, or at something that is considered a set in some set theory.  He goes on to argue that rather than putting the object first, in order to capture the idea of what a collection is, it is important to look at the relationship holding the objects together rather than first considering the objects, then the relations between them in a mathematical structure.  For instance, it is preferable and more meaningful to consider the orbits of an action on a set, considering the relationship induced by the action, rather than first focusing on the underlying object and worrying whether or not it is a set.  It is claimed that this relation-object approach as opposed to the object-relation approach is a more efficient and in a way more legitimate mathematical avenue, as for most mathematicians after all it is the relation rather than the underlying object itself that is considered more interesting and useful.  In short, Paul Taylor considers the nuts and bolts foundational approach to be an obstacle to the understanding of mathematics\ \cite{Tay}:

\begin{quote}
Foundations have acquired a bad name amongst mathematicians, because of the reductionist claim analogous to saying that the atomic chemistry of carbon, hydrogen, oxygen and nitrogen is enough to understand biology. Worse than this: whereas these elements are known with no question to be fundamental to life, the membership relation and the Sheffer stroke have no similar status in mathematics. 
\end{quote}

Just as simplifying biology to the point of chemistry is missing a larger picture, dealing with mathematical objects as simply a completed amalgamation of elements obfuscates the larger (and more important picture).  Namely, that a collection is determined by a predicate or relationship that should have priority.  To illustrate, one could consider $\mathbb{N}$ as the collection of anything that considers itself to be a "card carrying" natural number, or the collection defined by the relationship that every element is zero, or is one more than a preceding natural number.  While these two objects are mathematically the same, Taylor argues that in general the second approach is more instructive and in general preferable.

Of course, not everyone holds skepticism towards humanity's ability to properly reason with the completed infinite.  A prominent individual that would disagree with Taylor's skepticism towards the completed infinite is the logician Harvey Friedman.  From Friedman's perspective, completed infinities are necessary for understanding and practicing modern mathematics, and it is indeed preferable to consider these mathematical infinities as completed philosophical infinities.  Since exploring completed infinities through the charting of the transfinite or large cardinal axioms has immediate implications for propositions that might not at first glance appear to be related to each other, it is desirable to treat these infinities as completed objects that exist legitimately.  Succinctly, since we consider the ideas referred to by propositions as existing in some way, if the existence of a cardinal that is equivalent to this proposition or at least is implied by the proposition, it makes sense to accede that this cardinal is existing in some way that is at least as legitimate as the proposition in question.  This relationship is witnessed by the equivalence between some large cardinal propositions and issues in topology and combinatorics.

Part of Friedman's program is to demonstrate results using strong axioms, such as large cardinal axioms or axioms that are independent of the rest of the axioms of normal mathematics.  He believes, as do others, that many unresolved problems can be resolved, or at least fenced off as undecidable, using the relative consistency results that arise from large cardinals.  It is necessary to consider them as complete infinities in order to believe in the models required for these results, meaning that cordoning off the mathematical infinite as potential is insufficient for what Friedman identifies as a growing trend in mathematics.  Specifically, the use of strong axiomatic assumptions to determine the location of certain statements within a hierarchy of cardinals, acting as signposts along a path of increasing alienation from finite mathematics.

\section{G{\"o}del's Platonism}

\hspace{4ex}Kurt G{\"o}del is one of the most famous mathematicians of the 20th century, and arguably the most important.  His famous delimitative results, the incompleteness theorems, can be paraphrased as stating that any sufficiently strong formal system cannot prove its own consistency without being inconsistent.  The incompleteness theorems also state that a powerful enough system cannot decide, given a random proposition, whether it is provable or not.  The incompleteness theorems were the final nail in the coffin for the ideals of the positivists and logicists that proceeded before him, along the lines of the Hilbert program that sought to provide a firm axiomatic basis for mathematics that could be used to formalize the proof of any theorem that one pleased.
Positivism is the belief that all propositions are decidable, and logicism is the opinion that all mathematics is reducible to logical operations.
These theorems were the culmination of decades of progress in this direction, perhaps originally sparked by the infamous Russell's paradox.  Quickly, Russell's paradox arises from the fact that any predicate can define a set:  Is the set of all sets that do not contain themselves contained in itself?  Both possibilities lead to contradiction.  This paradox, which arises from disposing of any limitations on the predicates used to define sets, showed that there were unresolved foundational problems in mathematics.  

G{\"o}del himself held very nuanced positions on the reality of mathematics, and is generally considered as a Platonist as far as the metaphysics of the subject goes.  Similarly to Cantor, G{\"o}del had strong religious beliefs, but he did not use these beliefs as a tool to defend his philosophical position towards mathematics in contrast.  In fact, G{\"o}del worked on a formal and axiomatic proof of the existence of God via an ontological argument similar to that of St. Anselm.  Wary of the controversy that it might ignite, G{\"o}del kept it a secret until he thought he was dying\ \cite{Gode}.

G{\"o}del's was both a realist and a rationalist, in addition to being what could be called a Platonist. 
A realist is someone who believes in some sort of eminent reality of mathematics, and his rationalism could be interpreted as a restricted form of logicism.
Succinctly, his rationalist tendencies refer to the fact that he believed that "...under the impression that after sufficient clarification of the concepts in question it will be possible to conduct these discussions with mathematical rigour and that the result will then be…that the Platonistic view is the only one tenable\cite{Gode}".  This is an optimistic viewpoint that could be interpreted as an outgrowth of logical positivism, the idea that any problem is solvable in an appropriate logical system.  In light of his constraining results, namely the incompleteness theorems, it would perhaps seem anachronistic that G{\"o}del has such an optimistic viewpoint.  However, it is important to keep in mind that G{\"o}del would put restrictions on what problems are worth consideration.  Perhaps it could be rephrased as saying that any appropriate problem is possible with enough clarification, disposing of problems such as a system referring to its own consistency within itself.  An example of a meaningful proposition according to G{\"o}del would be the continuum hypothesis.  With enough clarification, with this term perhaps referring to explicitly choosing a mathematical model, this problem could be solved as it is in the constructible universe.  

His realist views can be summarized as the belief that mathematics is a reflection of a higher reality, not the reality that we sense perhaps with our eyes or ears, but a supposed reality that causes these sensory experiences.  The order and beauty found in mathematics is not due to some innate property of human reasoning, but the fact that this independent reality that causes our experiences is itself perfect and ordered.  This view is perhaps rooted in the idealism of Berkely 

The incompletenes theorems are often touted as arguments for a Platonistic conception of the reality of mathematics.  Namely, the incompleteness theorems could be used as a defense of the idea that mathematical objects exist independently of human existence, always have, and if one is a materialist could be taken further to the point of view stating that one of the views of physical reality, such as string theory, exist in a real way that determines material phenomena.  From a foundational point of view, this could be interpreted as a metaphysical claim that all things are numbers, collections, sets, functions, etc\ \cite{Port}.

How are G{\"o}del's incompleteness theorems used to bolster the apologetics of the Platonists?  Aleksander Mikovic argues that these theorems fall in line with expectations of science, by appealing that there must be undecidable phenonema in nature so that the incompleteness theorems must be a part of any sufficiently sophisticated physical theory.  In addition, he believes that some undecidable phenonema in physics arise from the fact that the passage of time is something that is not easily encapsulated in mathematics.  Mikovic's position could be summarized as stating that Platonism is a preferable viewpoint towards mathematics and physics, because of the fact that the incompleteness theorems capture abstract behaviors that are considered intuitive and fitting with our perception of reality, at least to Mikovic.  Here, mathematics has an abstract witness and a concrete witness in reality, with the perhaps boundless complexity of the universe providing propositions that are undecidable, and our finite perception of existence precluding the possibility of directly accessing the consistency of a sufficient physical theory.  

\begin{quote}

Since in the spirit of science is the belief that the natural laws can be described by mathematics,
then Platonism is a natural metaphysics for science. The natural laws must be mathematical,
because by definition, a natural law must be expressible by a finite string of symbols obeying the
laws of logic. One can also have a natural law within a materialistic metaphysics, but its meaning is
completely different from the meaning of a natural law in a platonic metaphysics. While a natural
law in a platonic metaphysics represents a timeless order in the motion of matter, in a materialistic
metaphysics a natural law is a temporary order which appears in the chaotic motion of matter at
random and lasts for a long time. Although a logical possibility, one can argue that this scenario is
highly implausible. Hence a platonic metaphysics is more plausible than a materialistic
metaphysics.

\end{quote}

\section{Formalism and mathematical utilitarianism}

\hspace{4ex}Formalism is an approach towards mathematics that emphasizes that syntactical properties of mathematics over any other sort of claim to reality.  From the formalist point of view, and variety of mathematics is legitimate and acceptable as long as the propositions are viewed as sequences of symbols that are derived methodically from base assumptions, at least while doing mathematics.  Some formalists are Platonists or finitists, but the key idea is that focusing on the logical and deductive aspects of mathematics is what is important, while passing over any other metaphysical claims in silence.  As long as mathematical rigor is preserved, then the formalist will acknowledge legitimacy as manipulation of symbols, but perhaps not give it any other kind of mathematical reality.  

Related to formalism, is the idea of mathematical utilitarianism.  Perhaps rooted in Cantor's appeal to mathematics being justified by its apparent consistency, mathematical utilitarianism can be encapsulated in the idea that the mathematics worth considering is mathematics that appears to be consistent, via relative consistency results and the existence of models in larger theories, and has the most applications, be it within itself, to other mathematical theories, or outside of mathematics altogether.  Joel Hamkins in particular applies this heuristic to the consideration of different varieties of set theory in his idea of infinitary utilitiarianism.  Specifically, the idea is that strong cardinal axioms are justified in their relation to other flavors of set theory, with those that are useful and producing the most interesting results as being most worthy of investigation\ \cite{Hamk}.  This pruning approach is especially desirable when the ideas often appear to be separated from normal human intuition and it is unclear which avenues to pursue when considered in an impartial vacuum.  

\section{Concluding remarks}

\hspace{4ex}While over the years there have been many disparate views towards the existence of mathematics, in the many ways this can be interpreted, it could be said that there are two fundamental groups that emerge.  In one corner there are those that choose their position by examining and considering the various metaphysical possibilities, particularly the epistemological status that mathematics would be left in, and choose the belief system that leads most closely to what they considered to be the appropriate position of mathematics in the hierarchy of human activities.  On the other hand, there are those that first start with notions about reality that are considered intuitive and attempt to inductively build up these ideas into a scheme on the existence of mathematics, perhaps leaving it in a place independent of human experience. 

This fundamental difference in motivation does not necessarily lead to incompatible conclusions.  An individual from one camp might reach the same position as their counterpart in the other, but for very different reasons.  Ultimately, everyone is choosing a position on the topic of mathematics, and I claim that this separation indicated in the previous paragraph is a major motivational difference that can be summarized as the distinction between placing mathematics or reality as it is sensed as the first and prioritized idea in what we should think about mathematics. Alternatively, this difference could be portrayed as the distinction between considering only mathematics when thinking of what position mathematics will hold, or trying to reach a verdict on the status of mathematics by considering it in a spectrum of philosophical ideas, ranging from real experience to non-mathematical abstract reasoning and intuition.

Foremost among those who consider the mathematics first and form their beliefs about mathematics based on what happens after considering different possibilities is early Cantor before his religious turn.  Early Cantor's requirement for mathematical existence was for a theory to be clear and at least have the veneer of consistency: as long as the theory seems like it makes sense, it is good to go and as legitimate as any other, regardless of the implications outside of mathematics.  Here, early Cantor is a Platonist and a mathematical realist, and thus can be considered to be in the same camp as G{\"o}del, at least as far as the metaphysics of mathematics goes.  

On the other hand, G{\"o}del reached his realist and Platonist conclusions for different reasons than Cantor, reasons that place him in the camp that puts reality, philosophy, and experience as foundations for how we think of mathematics.  Here I am emphasizing the difference between a foundation for how we think of mathematics and a foundation for mathematics itself, such as in some set, type, or category theory.  For instance, rather than referring to the apparent consistency of mathematics as a justification for its real existence as Cantor does, G{\"o}del mirrors his ontological argument for God as an ontological argument for the perfect existence of mathematics as an objective and crystalline reality: "We may expect that the conceptual world is perfect, and, furthermore, that objective reality is beautiful, good, and perfect."  So, Cantor and G{\"o}del reach the same destination, with the former arriving via an almost utilitarian argument, and the other by an ontological appeal to the purported reality of a perfect abstract space.

On the other hand, there are those that reach different metaphysical conclusions on mathematics, while tapping in to the same driving motivation as the other.  For instance, both late religious Cantor and Kronecker used arguments outside of mathematics to bolster the position that mathematics exists as a real and indifferent space, but as is mentioned early had vitriolic disagreements over the reality of specific parts of what is considered today to be mathematical canon.  Ultimately though, despite all of these differences, be it in motivation or conclusion, all of the people that are mentioned above are individuals who have contributed in some way to the progression of mathematics, throwing in their contributions, small or large, to the project of mathematics.  It is important to remember this, and perhaps it is one of the great accomplishments of mathematics outside of itself that individuals with incompatible viewpoints are able to fruitfully collaborate towards the higher goal of expanding the breadth of human knowledge.

\end{document}